\documentclass{article}
\usepackage[english]{babel}
\usepackage{amssymb}
\usepackage{pstricks}
\usepackage{epic,eepic,graphicx}
\usepackage{amsmath}

\begin{document}

\title{{\huge On the history of ring geometry \\(with a thematical overview of literature)}}

\author{{\Large Dirk Keppens} \bigskip \\Faculty of Engineering Technology, KU Leuven\\Gebr. Desmetstraat 1 \\B-9000 Ghent 
BELGIUM \\\textbf{\texttt{dirk.keppens@kuleuven.be}}}

%\email{}

%\subjclass{Primary 51C05; Secondary 00A15, 01A60, 01A61, 51-03, 51A05, 51B05, 51E26, 51H10, 51F15}

%\keywords{Ring geometry, projective ring plane, Hjelmslev geometry, Klingenberg geometry, Barbilian plane, neighbor relation, bibliography}

\date{}

%\dedicatory{\sl Dedicated to Willem Mielants on the occasion of his 75th birthday,\\ in friendship and honour.} 

\maketitle

\begin{abstract}
In this survey paper we give an historical and at the same time thematical overview of the development of ``ring geometry" from its origin to the current state of the art. A comprehensive up-to-date list of literature is added with articles that treat ring geometry within the scope of incidence geometry. \\

In questo documento di ricerca forniamo una panoramica storica e allo stesso tempo tematica dello sviluppo della ``geometria sopra un anello" dalla sua origine allo stato attuale. \`E aggiunto una lista di letteratura aggiornata completa di articoli che trattano la geometria degli anelli nel contesto della geometria dell'incidenza.\\

In diesem Forschungsartikel geben wir einen historischen und gleich\-zeitig thematischen \"Uberblick \"uber die Entwicklung der ``Ringgeometrie" von ihrem Ursprung bis zum aktuellen Stand der Technik, mit einer Liste aktualisierter Literatur einschlie\ss lich Artikeln zur Ringgeometrie im Kontext der Inzidenzgeometrie.\\

Dans ce document de recherche, nous fournissons un aper\c{c}u historique et \`a la fois th\'ematique du d\'eveloppement de la ``g\'eom\'etrie sur un anneau", de son origine \`a l'\'etat actuel des connaissances. Nous ajoutons une liste de la litt\'erature actualis\'ee comprenant des articles traitant la g\'eom\'etrie sur un anneau dans le contexte de la g\'eom\'etrie de l'incidence.
\end{abstract}

\bigskip


\noindent


\section{Introduction}
The current version of the Mathematics Subject Classification (MSC), a classification scheme used by the two major mathematical reviewing data\-bases, Mathematical Reviews (MathSciNet) and Zentralblatt MATH (zbMATH), provides code 51C05 for all research papers dealing with {\sl Ring Geometry}. This rather small category contains articles about geometries that are not only provided with an incidence relation but also with a neighbor relation (or its negation, a non--neighbor or remoteness relation). Included are all geometries obtained from rings that are not division rings. Ring geometry, i.e.~the \mbox{theory} of geometries equipped with a neighbor relation in general and of geometries over rings in particular, is a rather young discipline. Its origin lies in the beginning of the 20th century and its importance has grown steadily. For that reason it has also got a full--fledged place as a chapter in the {\sl ``Handbook of Incidence Geometry"} \cite{VeldL}. In the past decades \mbox{several} mathe\-maticians have contributed to ring geometry. Newly discovered connections with \mbox{coding} theory, with the theory of Tits--buildings and with quantum information theory, have opened new horizons. In this survey paper we present an historical and thematical overview with attention to many aspects. An out-of-date list with articles on the subject, up to 1989, was composed by {\sc T\"orner} and {\sc Veldkamp} in \cite{TV} as an addition and completion of two even older lists \cite{ADDT} and \cite{JungL} written by respectively {\sc Artmann, Drake} {\sl et al.} and {\sc Jungnickel}. Up to now such list of literature is not available for the period after 1990, except for a survey paper by the author \cite{KepL}, dealing exclusively with plane projective geometries over finite rings. In the present work we fill that gap and we add a new updated list of existing literature, ordered thematically (and containing the relevant material from the \mbox{preceding} lists). Articles on algebraic geometry and on differential geometry over rings are not included. We also do not aim for completeness when it concerns metric aspects, neither for geometries on modular lattices nor for geometric algebra over rings. We think that the bibliography might be useful for researchers who want to attack the future challenges of ring geometry.  

\section{The first traces of ring geometry: dual numbers and Johannes Hjelmslev}
The first traces of ring geometry date back to the beginning of the twentieth century.
The Danish geometer Johannes Trolle {\sc Hjelmslev} (1873--1950) who was born as Johannes Petersen but who changed his name in 1904, graduated in mathematics from the University of Copenhagen and received his PhD degree in 1897.  Hjelmslev can be viewed as one of the early founders of ring geometry. In a series of four lectures, held at the University of Hamburg in July 1922 and published one year later in  the Abhandlungen aus dem Mathematischen Se\-minar \cite{Hjelmslev2}, he presented an axiomatic framework for geometry that better reflected the properties observed in the real world. The basic observation made by Hjelmslev was that, if one draws lines ``close" to each other (meaning that the sharp angle they define is very small), then it is hard to identify the intersection point, and it looks as if the lines have a little segment in common. Dually, if two points are close to each other, they belong to a line segment that can be part of several joining lines. 
Hjelmslev called this {\sl Die nat\"urliche Geometrie}, the ``natural geometry".
In fact, Hjelmslev put forward this idea already some years earlier, in 1916, in an article \cite{Hjelmslev1} on what he called {\sl Die Geometrie der Wirklichkeit}, the ``geometry of reality".

In order to obtain a model for his geometry, Hjelmslev made use of the ring of dual numbers of the form $a+b\varepsilon$ with $a$ and $b$ both real  and $\varepsilon^2=0$. \\ Dual numbers were already well--known before Hjelmslev. William {\sc Clifford} defined them for the first time in 1873 \cite{Clifford} and they were used as a convenient tool in mechanics by Aleksander {\sc Kotelnikov} \cite{Kotelnikov} and by Eduard {\sc Study} in his famous work {\sl ``Geometrie der Dynamen"} \cite{Study}.\\
According to Benz \cite{BenzF} the first traces of ring geometry can be observed already in Study's work and in that of some of his contemporaries like Josef {\sc Gr\"unwald},  Pilo {\sc Predella} and Corrado {\sc Segre}. In their papers dual numbers are treated from a geometrical viewpoint, by connecting the geometry of oriented lines (spears) in the real euclidean space with a spherical geometry over the ring of dual numbers \cite{Grunwald,Pred,SegreC}.

\section{The pioneers of ring geometry: Barbilian and Klingenberg}

Twenty years after Hjelmslev, Dan {\sc Barbilian} (1895--1961), besides a mathe\-matician at the University of Bucharest, also one of the greatest Romanian poets (with pseudonym Ion Barbu), took up the line again. In two papers \cite{Barbi1, Barbi2} (extended versions of a lecture hold in Baden--Baden in 1940)  he gave the first axiomatic foundation for plane projective ring geometry. He started by investigating the conditions which must be imposed upon an arbitrary associative ring in order that the corresponding geo\-metry may have ``useful"  properties. It turns out that the rings must have a unit element and that all left singular elements must be two--sided singular or equivalently  $ab = 1$ implies $ba = 1$. Barbilian called these rings ``Zweiseitig singul\"are Ringe" or $Z$--rings. Today such rings are also known as Dedekind--finite rings. They include of course all commutative rings but also all finite rings (even non--commutative) and many other classes of rings (e.g.~matrix rings over a field). Starting with $C$--rings (``Kategorische Ringe"), being $Z$--rings with some additional property, Barbilian defined a kind of projective plane. Conversely, he formulated a set of axioms for a plane geometry with incidence and non--neighborship (``klare Lage") and proved that this leads to a $C$--ring. 
Barbilian's work showed some shortcomings, but nevertheless it was of great importance for the development of geometry over arbitrary rings as we will  discuss further (see also section 11). 

Another monument among the ring geometers was the German mathe\-matician Wilhelm {\sc Klingenberg} (1924--2010). He is best known for his work on diffe\-rential geometry, but in a series of papers \cite{Kling1,Kling2,Kling3,Kling4}, \mbox{published} in the mid--fifties, he also laid the foundation for affine, metric (euclidean) and projective geometries with neighbor relation (``mit Nachbarelementen"). His central idea is the existence of a natural map from the geometry under consideration onto an ``underlying" ordinary geometry, a consequence of the assumption that the neighbor relation is transitive. Therefore Klingenberg called such geometries ``Geometrien mit Homomorphismus". He not only defined them axiomatically, but he also constructed explicit examples using rings. The assumption of a transitive neighbor relation is interlaced with the fact that the rings must be local ones. A not necessarily commutative local ring $R$ is a ring with a unique maximal left (or equivalently a unique maximal right) ideal  $R_0$. The quotient ring $R/R_0$ is a (skew)field, coordina\-tizing an ordinary geometry which is the natural epimorphic image of the ring geometry over $R$. 

The work of Klingenberg was the source of two mainstreams in ring geo\-metry in the decades afterwards: {\sl Klingenberg geo\-metry} and {\sl Hjelmslev \mbox{geometry}}. Geometries with a transitive neighbor relation or equivalently with a canonical map onto an ordinary geometry are now named Klingenberg geo\-metries. Their epimorphic image is obtained by considering the equivalence classes of neighboring elements as ``thick" elements. Non--neighboring elements behave like distinct elements in the epimorphic image (e.g.~two non--neighboring points are incident with exactly one line in a projective Klingenberg plane since two distinct points in its epimorphic projective plane are incident with a unique line). No further assumptions are made about the neighboring elements.\\ Klingenberg geometries comprise the smaller, but more intensively studied class of Hjelmslev geometries for which additional axioms must hold when it comes to neighboring ele\-ments (e.g.~two neighboring points can be connected by at least two distinct lines in a projective Hjelmslev plane). This idea relies on the natural geo\-metry of Hjelmslev (see section 2). In the case of a Klingenberg plane over a local ring $R$, the Hjelmslev condition means that the ring must be a two--sided chain ring (a ring is a right chain ring if for any $a$ and $b$ in $R$ either $a \in bR$ or $b \in aR$, similar for left chain ring) with the additional property that its maximal ideal contains only zero divisors. Such rings are called PH--rings (projective Hjelmslev rings or $H$--rings). Finite chain rings are always $H$--rings and, equivalently, local principal ideal rings. Klingenberg himself considered both Klingenberg planes (over local rings) and Hjelmslev planes (over Hjelmslev rings) in \cite{Kling3}. 

To be complete, we must mention also two isolated, rather unnoticed, contributions to early geometry over rings. J.W. {\sc Archbold}  \cite{Archbold} wrote on projective geometry over  group algebras $K[G]$ and worked out a small finite example \mbox{taking} $K=$ GF(2) and $G$ a group of order 2,
while Rimhak {\sc Ree} in \cite{Ree} considered projective spaces as the modular lattices $L(R,M)$ of all sub--modules of a module $M$ over a ring $R$, in particular over a full matrix ring.

\section{The Belgian contribution to early ring \mbox{geometry:} projective geometry over rings of matrices}

The ideas of Barbilian concerning projective geometries over rings instead of fields, got some followers in the sixties and early seventies. Among them we find a lot of Belgian contributors who studied projective geometries over rings of matrices. Since such rings are not local, the geometries are not Klingenberg nor Hjelmslev geometries.

Julien {\sc Depunt} studied in \cite{Bilo} the projective line and in \cite{Depunt1,Depunt2} projective planes over the ring of ternions. This ring contains the elements $a_1\varepsilon_1+a_2\varepsilon_2+a_3\varepsilon_3$ with $a_i \in \mathbb{C}$ and  $\varepsilon_i^2=\varepsilon_i$ ($i=1,2$),\, $\varepsilon_1\varepsilon_3=\varepsilon_3=\varepsilon_3\varepsilon_2$ and all other products equal to $0$. The ring of ternions is isomorphic to the ring of upper triangular matrices with ele\-ments in the complex field. The main result of Depunt concerns the embedding of the plane over the ternions into the 5--dimensional complex projective space.

Inspired by this work,  Cl\'ery {\sc Vanhelleputte} could generalize the results for planes over the full matrix ring of $2 \times 2$--matrices with elements in an arbitrary commutative field. His doctoral thesis was published in \cite{Vanhelleputte} and fitted in the research program under Julien {\sc Bilo} (1914--2006) who was head of the Geometry Department at the University of Ghent at that time.

A few years later, in 1969, Joseph {\sc Thas} wrote his PhD thesis under the supervision of Bilo about the projective line over the ring of $3 \times 3$--matrices with elements in an algebraically closed field \cite{Thas2}, see also \cite{Thas1}. Soon afterwards Thas published \cite{Thas3,Thas4} on projective ring geometry over matrix rings, in particular over finite rings $M_n(GF(q))$ of $n \times n$--matrices with entries in a Galois field of order $q$. Especially \cite{Thas3} was important as it contains some concepts (ovals, arcs, caps, etc.) and combinatorial theorems which extend known results from classical Galois geometry over finite fields and which appear for the first time for ring geometries. Later Thas grew out to an authority in finite geometry, but his interest shifted from finite ring geometry to Galois geometry and other \mbox{topics}, in particular finite generalized quadrangles.\\ Thas was a colleague of Willem {\sc Mielants} during many years. It was he who guided Hendrik {\sc Van Maldeghem} (see section 14) and the author towards ring geometry as PhD super\-visor for both of us.

Paul {\sc De Winne} was the last student of Bilo studying ring geometry, see \cite{DeWinne}. His doctoral thesis was published in \cite{DeWinne2}.

Around the same time, Xavier {\sc Hubaut}, another Belgian mathematician from the Universit\'e Libre de Bruxelles, introduced the projective line over a general associative algebra and investigated the structure of its group of projectivities in \cite{Hubaut1,Hubaut2}. His colleague Franz {\sc Bingen} could generalize the results to $n$--dimensional projective spaces over semi--primary rings in \cite{Bingen}. 

Mind that the geometries over matrix rings studied by all these people, may not be confused with the ``geometries of matrices" initiated by the Chinese mathematician Loo--Keng {\sc Hua} in the mid--forties \cite{Hua}. Hua's geometries of matrices are not ring geometries in the narrow sense, despite the suggestive name (although a connection is possible with projective lines over rings, see section 13). 

\section{The foundations of plane affine ring geometry: from Benz to Leissner and beyond}

Barbilian, as well as all members of the Belgian School, dealt exclusively with {\sl projective} ring geometry. The first traces of {\sl affine} geometries over rings can be found in a paper from 1942 by Cornelius {\sc Everett}, concerning affine planes over rings without zero divisors \cite{Everett}.\\  Walter {\sc Benz} (1931--2017), a renowned German geometer, especially for his work on circle geometries, considered in \cite{Benz1,Benz2} affine planes over a special kind of commutative rings and in \cite{Benz3} plane affine (and metric--affine) geometries over arbitrary rings with unit. A big part of this interesting paper discusses the relation between algebraic properties of the ring and geometric properties of the plane.

Peter Ashley {\sc Lawrence} obtained a PhD under the supervision of Benz. His thesis {\sl ``Affine mappings in the geometries of algebras"} considers ring geometries associated with modules and algebras over a commutative ring and was published in \cite{Lawrence}.
Two other publications \cite{Arnold2,Arnold1} by Hans--Joachim {\sc Arnold} can be classified within this section. They contain a simultaneous generalization of affine geometry over rings and generalized affine spaces in the sense of Emanuel {\sc Sperner} \cite{Sperner}. Arnold starts with the axiomatic defi\-nition of an affine--line geometry which can be coordinatized by a vectorial groupoid and he adds some extra axioms to turn it into a module over a unitary ring. Also \cite{Kaerlein} and \cite{Schleicher}  both fit into this approach.

Affine planes over local rings, today also known as desarguesian affine Klingenberg planes, appear as particular examples in the papers of Benz, but this was not for the first time. {\sc Klingenberg} already considered them on the side of projective ring planes and in the more general context of ``Affine Ebenen mit Nachbarelementen" in the articles \cite{Kling1,Kling2,Kling3,Kling4}.\\ A reasonable part of the research concerning  affine ring geometry focusses on affine Hjelmslev planes (AH--planes). 
An important paper \cite{Luneburg} in that respect was written by Heinz {\sc L\"uneburg} (1935--2009) in 1962. It contains the main results of his doctoral thesis entitled {\sl ``Affine Hjelmslev--Ebenen mit transitiver Translationsgruppe"}, dealing with translation AH--planes as generalizations of ordinary translation planes. Werner {\sc Seier} investigated in more detail desarguesian and translation AH--planes in \cite{SeierA1,SeierA2,SeierA3,SeierA4,SeierA5}. The central theme in his work is a characterization of desarguesian AH--planes by an affine variant of Desargues' theorem and the transitivity of the set of translations in the plane. Some other papers on affine Hjelmslev planes are not mentioned here and we postpone them to sections 7 and 8 where they will be discussed in connection with projective Hjelmslev planes.

Without minimizing the importance of the papers cited above, we can state that the real breakthrough of affine ring geometry in full generality was not achieved yet at this point. Therefore one has to wait for Werner {\sc Leissner}, a student of Benz. Leissner wrote two papers in the mid--seventies on, what he called, ``affine Barbilian planes" (today the name ``Leissner planes" would be more convenient). In \cite{Leissner1} such planes are defined axioma\-tically as affine structures in which two parallelodromy axioms hold (suitable substitutes for the affine Desargues theorem) and they are coordinatized by an arbitrary $Z$--ring $R$ together with a special subset $B$ of $R \times R$ (a Barbilian domain). In \cite{Leissner2} the converse is proved: any affine ring plane over a $Z$--ring is an affine Barbilian plane. Leissner considered even more \mbox{general} affine Barbilian structures in \cite{Leissner3}.\\ The papers of Leissner were a source of inspiration for further generalizations by several people. Victoria {\sc Groze} and Angela {\sc Vasiu} con\-sidered affine structures over arbitrary rings in \cite{GV}, Francisc {\sc Rad\'o} defined affine Barbilian structures in \cite{Rado1} by weakening Leissner's axioms, {\sc Armentrout} {\sl et al.} investigated in \cite{AHM} generalized affine planes which can be coordinatized by near-rings and {\sc Pickert} studied tactical configurations over quasirings in \cite{Pickert1}.

Other contributors to the theory of plane affine geometries with neighbor relation and parallellism were Gernot {\sc Dorn} \cite{Dorn}, Kvetoslav {\sc Burian} \cite{Burian1,Burian2,Burian3}, Frantisek {\sc Machala} \cite{MachalaA1,MachalaA2,MachalaA3}, Stefan {\sc Schmidt} and Ralph {\sc Steinitz} \cite{SS1}, Franco {\sc Eugeni} \cite{EuGa,EuMa} and Angela {\sc Vasiu} \cite{Vasiu1,Vasiu2,Vasiu3}.
The Barbilian domains defined by Leissner also formed a study object in themselves. Several authors have contributed to this, see \cite{Benz4, Lantz, Leissner4, Leissner5,LorimerZ}.
For a recent and fairly complete overview of affine planes over {\sl finite} rings we refer to the survey paper \cite{KeppensA} of the author.
We also mention here two papers \cite{CelikH1,CelikH2} by Basra {\sc \c{C}elik} on finite hyperbolic Klingenberg planes as they are not far away from finite affine Klingenberg planes.

\section{Metric geometry over rings and the school of Bachmann}

The goal of this section is just to give a glimpse of metric ring geometry. We do not aim for completeness when it concerns metric aspects, because they are only indirectly related to incidence geometry. The literature about this subject is extensive and we have selected only a few representative papers in which much more references can be found. The basis of ``modern" metric geometry is provided by the work of several German mathematicians, including Bachmann, Lingenberg and Schr\"oder.\\
The pivotal figure is Friedrich {\sc Bachmann} (1909--1982). In the second edition of his famous book {\sl ``Aufbau der Geometrie aus dem Spiegelungsbegriff"} \cite{Bachmann1} he considers plane metric geometries in which points, lines, incidence and orthogonality are defined by means of a group of reflections. 
A group $G$ with a subset $S$ of involutory elements which are invariant under inner automorphisms of $G$ and which generate $G$, determines a group plane $E=E(G,S)$ as follows: the elements of $S$ are the lines of $E$ and the involutory elements of $S^2$ are the points. Two lines that commute are perpendicular. A point and a line that commute are incident.  If the group plane $E$ satisfies the conditions (1) a point and a line determine a unique perpendicular, and (2) the product of three lines that are either concurrent or have a common perpendicular is again a line, then the pair $(G,S)$ is called a Hjelmslev group and the associated group plane $E$ is called a {\sl metric (non--elliptic) Hjelmslev plane}. \\If the uniqueness of the line incident with two distinct points is not required, then metric Hjelmslev planes stripped of their metric structure, reduce to incidence Hjelmslev planes.  In section 2 we have already discussed the role of Hjelmslev as founder of (incidence) ring geometry. In some later work, the {\sl ``Algemeine Kongruenzlehre"} \cite{HjelmslevM}, Hjelmslev used special transformations (reflections) to define ortho\-go\-nality. Hence, the rudiments of metric ring geometry also go back to Hjelmslev and have been worked out in more detail by Bachmann. There is also a minor contribution by Klingenberg who treats metric aspects in ring geometries in \cite{Kling2}.

Rolf {\sc Lingenberg} (1929--1978), a student of Bachmann, continues his work in \cite{Lingenberg}. The main result is an algebraic characterization of classical metric group planes $E(G,S)$ with an additional axiom, as planes associated with a \mbox{metric} vectorspace $(V,q)$ with $q$ a quadratic form.\\Eberhard {\sc Schr\"oder} in \cite{Schroder1, Schroder2} considers metric planes of different kind starting from a pappian affine plane over a field and using two\-dimensional algebras for the introduction of the metric notions of angle, distance and orthogonality. His work is strongly related with circle geometries studied intensively by Benz (see section 13). A good overview of classical \mbox{metric} geometry (but not including  ring geometries) can be found in Chapter 17 written by Schr\"oder in the {\sl Handbook of Incidence Geometry} \cite{Schroder3}.

There are numerous contributions by Bachmann and his school to the theory of Hjelmslev groups, connecting geometric and algebraic properties of groups and planes.
From Bachmann himself we mention the important additions \cite{Bachmann3, Bachmann2, Bachmann4} to his book \cite{Bachmann1}.
In the seventies and eighties his successors, most of them from the University of Kiel (Germany), have extended the knowledge. 
Finite Hjelmslev groups were characterized by Rolf \mbox{{\sc St\"olting}} in his PhD thesis {\sl ``Endliche Hjelmslev-Gruppen und Erweiterungen von Hjelmslev-Gruppen"}, also published in \cite{Stolting1}. In \cite{Stolting2} Hjelmslev groups are constructed from a module $M$ over a commutative ring $R$ endowed with a bilinear form. \\
Edzard {\sc Salow} introduced singular Hjelmslev groups (in which the \mbox{product} of three points is always a point) in his doctoral thesis {\sl ``Beitr\"age zur Theorie der Hjelmslev-Gruppen: Homomorphismen und Singul\"are Hjelmslev-Gruppen"} published in \cite{Salow1}. The main result is the construction of a coordinate ring $R$ for the group plane of a singular Hjelmslev group, proving that these \mbox{metric} planes are indeed ring geometries. 
The process of algebraization of metric Hjelmslev planes is investigated also in \cite{Salow2}. It is proved there that a Hjelmslev group  with some additional axioms can be embedded in the orthogonal group of a metric module $(R^3,f)$ with $R$ a commutative ring with unit for which $2$ and any non zero--divisor is invertible and with $f$ a sym\-metric bilinear form of the free module $R^3$.\\
In \cite{Salow3} Salow studies another class of metric ring planes using a commutative algebra over a ring and the concept of an angle. In an early paper of Benz \cite{Benz1} the \mbox{metric} notion of angle also played an important role and in \cite{Benz3} he devotes a paragraph to metric geometry, using an elliptic form over an arbitrary commutative ring.\\
Similar work can be found in \cite{Nolte1, Nolte2} in which Wolfgang {\sc Nolte} proves that a class of metric planes $E(G,S)$ with additional axioms can be embedded into a projective Hjelmslev plane over a local ring and that $G$ is isomorphic to a subgroup of an orthogonal group. Other results in this direction were obtained by Frieder {\sc Kn\"uppel}, in \cite{Knuppel2,Knuppel3,Knuppel4,Knuppel5,Knuppel6}, Gerald {\sc Fischbach} \cite{Fischbach}, Michael {\sc Kunze} \cite{Kunze} and Rolf and Horst {\sc Struve} \cite{Struve1,Struve2,Struve3,Struve4}.\\
The influence of Bachmann is also clear from the large amount of doctoral theses produced in Germany on Hjelmslev groups and \mbox{metric} geometry, e.g.~R. Schnabel (1974), M. Kunze (1975), H. Struve (1979), R. Struve (1979), M. Gegenwart (1987), W. Vonck (1988) and A. Bach (1998).

\section{The florescence of Hjelmslev geometry in the era of Drake, Artmann and T\"orner}

The first period of florescence of ring geometry (especially Hjelmslev geo\-metry) regarded as incidence geometry, started in the late sixties and reached his culmination point in the seventies. This is reflected in a large number of publications. Two mathe\-maticians who were very productive in this area and left their mark, were David Allyn {\sc Drake} (1937--) from the University of Florida, Gainesville (USA) and Benno {\sc Artmann} (1933--2010) from the University of Giessen (Germany).

Drake obtained his PhD in 1967 with a thesis entitled {\sl ``Neighborhood collineations and affine projective extensions of Hjelmslev planes"} under the supervision of Erwin Kleinfeld. {\sc Kleinfeld} himself, an authority in algebra, with a lot of publications on non--associative alternative rings, published only one, though interesting paper \cite{Kleinfeld} about ring geometry. It was the first publi\-cation concerning {\sl finite} Hjelmslev planes.
Among other things he introduced a two parameter set $(s,t)$ of non--zero integers such that for each flag $(P,\ell)$ in a finite projective Hjelmslev plane, there are exactly $t$ points on $\ell$ neighboring with $P$ and exactly $s$ points on $\ell$ not neighboring with $P$. It was proved that $s \leq t^2$ or $t=1$. If $s=t^2$, the plane is called uniform. In that case all point neighborhoods have the structure of ordinary affine planes.  Robert {\sc Craig} proved in \cite{Craig} that any finite projective plane can be extended to a uniform projective Hjelmslev plane. \\ 
The notion of uniformity (and its generalization to $n$--uniformity) has played a crucial role in the work of Drake. It was also related to another issue: the extension of an affine Hjelmslev plane to a projective Hjelmslev plane. An ordinary affine plane can always be extended to a projective plane, but for Hjelmslev planes the situation is much more complicated. Drake has written several papers about this problem in the period from 1968 to 1975. In \cite{Drake1} he proves that any uniform affine Hjelmslev plane has at least one (uniform) projective extension and he gives an example of a desarguesian uniform affine Hjelmslev plane with a non--desarguesian projective extension. In \cite{Drake2} uniformity is generalized to $n$--uniformity inductively (a PH--plane is $n$--uniform if the point neighborhoods are $(n-1)$-uniform AH--planes). Strongly $n$--uniform planes are characterized by a local property, which leads to the theorem: an $n$--uniform PH--plane is strongly $n$--uniform if and only if its dual is $n$--uniform. Drake also proved that every finite desarguesian PH--plane is strongly $n$--uniform. \\
A further study of $n$--uniform Hjelmslev planes (projective and affine) was made in \cite{Drake3,Drake4,Drake5,Drake6,Drake7,Drake8,Drake14} where even more general affine geometries with neighbor relation  appear. Drake could also prove that there do \mbox{exist} affine Hjelmslev planes which cannot be extended to projective Hjelmslev planes. 

Artmann was a contemporary of Drake and he wrote his doctoral thesis {\sl ``Automorphismen und Koordinaten bei ebenen Verb\"anden"} in 1965 under the supervision of G\"unther Pickert.
In his early work on Hjelmslev geo\-metry we can observe a strong relation with the theory of modular lattices. In \cite{Art1} he gives a sufficient condition for a modular lattice to define a projective Hjelmslev plane and in \cite{Art3} he proves that any uniform PH--plane can be derived from a modular lattice.\\ 
Like Drake, Artmann also studies refinements of the neighbor relation (``verfeinerten Nachbarschaftsrelationen") and the affine--projective extension question. In \cite{Art5} he proves that a uniform affine Hjelmslev plane can be extended to at least two non--isomorphic projective Hjelmslev planes. A new concept introduced by him in \cite{Art2} is that of a projective Hjelmslev plane of level $n$ (``$n$--ter Stufe") based on the refinement of  the neighbor relation. This definition was extended by Drake to the affine case in \cite{Drake7}.  Artmann proves that desarguesian PH--planes over a Hjelmslev ring $R$ are of level $n$ if and only if the maximal ideal of $R$ is nilpotent of index $n$, see \cite{Art6,Art7}. 

Another theorem proved by Artmann \cite{Art8} states that for any projective plane $\mathcal{P}$ and any integer $n>0$ there exists a PH--plane of level $n$ with $\mathcal{P}$ as epimorphic image. Moreover, given a sequence of PH--planes $\ldots \rightarrow H_i \rightarrow H_{i-1} \rightarrow \ldots \rightarrow H_1=\mathcal{P}$ with $H_i$ of level $i$, the inverse limit is a projective plane.
Arno {\sc Cronheim} constructed in \cite{Cronheim1} in a purely algebraic way, using formal power series over a cartesian group, a chain of Hjelmslev planes whose inverse limit is a projective plane.  Cronheim also obtained a complete classification of all finite uniform desarguesian projective Hjelmslev planes in \cite{Cronheim2}. They are either planes over a ring of twisted dual numbers over GF($q$) (a non--commutative generalization of the classical dual numbers) or over a truncated Witt ring $W_2(q)$ of length 2. 

Both Drake and Artmann had a great influence on the mathematical research in the domain of Hjelmslev geometry (even when Artmann's \mbox{interest} shifted to other subjects after a few years). One of Drake's students was Phyrne {\sc Bacon}. She wrote both her Master's thesis {\sl ``On Hjelmslev planes with small invariants"} and her PhD thesis {\sl ``Coordinatized Hjelmslev planes"} on Hjelmslev geometry, resulting in two papers \cite{Bacon1} and \cite{Bacon2}. She proved that a finite Hjelmslev plane is strongly $n$--uniform if and only if it is of level $n$ which unifies the two notions introduced by Drake and Artmann respectively. Later, Bacon's attention shifted to the more general Klingenberg geometries (see section 9). \\
Artmann was the supervisor of Manfred {\sc Dugas} and G\"unther {\sc T\"orner} who both made important contributions to Hjelmslev geometry. The PhD thesis of Dugas {\sl ``Charakterisierungen endlicher desarguescher uniformer Hjelmslev-Ebenen"} contains many new ideas, including a coordinatization method (see section 9). 
In \cite{Dugas2} Dugas proves that a finite translation AH--plane can be extended to a PH--plane if it is or a desarguesian PH--plane or an ordinary translation projective plane, while in \cite{Dugas5} he gives a necessary and sufficient condition for a projective Hjelmslev plane to be derivable from a lattice. In particular the PH--planes of level $n$ are always lattice--derivable.

T\"orner wrote his Master's thesis on {\sl ``Hjelmslev--Ringe und die Geometrie der Nachbarschaftsbereiche in der zugeh\"origen Ebenen"} and obtained his PhD with {\sl ``Eine Klassifizierung von Hjelmslev-Ringen und Hjelmslev-Ebenen"} in the same year as Dugas, under the supervision of Artmann and Pickert. His research in the domain of ring geometry focusses on two main themes: the structure of (finite) Hjelmslev planes and the ideal structure of chain rings. Among his publications we mention here \cite{Torner1} which contains the main results from his thesis: a classification of PH--planes based on congruence relations. He proves that the set of all congruence relations of a finite PH--plane is linearly ordered under inclusion and consequently, the canonical epimorphism onto the associated projective plane admits an essentially unique factorization into indecomposable epimorphisms. The plane  $\mathcal{H}$ is of ``type $n$" or ``height $n$" if the canonical epimorphism $\varphi$ from $\mathcal{H}$ onto the projective plane $\overline{\mathcal{H}}$ has a maximal factorization $\mathcal{H}=\mathcal{H}_n \rightarrow \mathcal{H}_{n-1} \rightarrow \ldots \rightarrow \mathcal{H}_1=\overline{\mathcal{H}}$. In \cite{Torner3,Torner4} T\"orner investigates the equivalence of finite $n$--uniform planes or planes of level $n$ as defined by Drake and Artmann to planes of type $n$. In \cite{Torner2} he proves that $n$--uniform projective Hjelmslev planes are strongly $n$--uniform. Some of the results were later extended to the infinite case in \cite{Torner9}. In \cite{Torner3} much attention goes also to affine Hjelmslev planes, in particular translation AH--planes over near--rings. In T\"orner's work homomorphisms play an important role as can also be seen from \cite{Torner5,Torner6}.

In the desarguesian case (Hjelmslev planes over a chain ring) the structure of the plane is intrinsically connected with the ideal structure of the ring. The structure of chain rings and valuation rings was investigated by T\"orner partly in collaboration with Hans--Heinrich {\sc Brungs}. One of their papers \cite{BrungsTorner2}  concerns the embedding of right chain rings into chain rings (related with the problem of embedding desarguesian affine Hjelmslev planes into projective ones), see also \cite{Mach1,Mach2,Skornjakov,Torner8}. \\
With a postdoctoral scholarship Artmann stayed for a short time at the McMaster University in Ontario (Canada). There he inspired J.W.~(Mike) {\sc Lorimer} who would later become one of the leading figures in topological Hjelmslev geometry (see section 10). In \cite{LoLa} Lorimer and Lane study desarguesian Hjelmslev planes. They prove that an affine Hjelmslev plane is desarguesian if and only if it can be coordinatized by an AH--ring and that not every desarguesian AH--plane can be extended to a desarguesian PH--plane. Morphisms between affine Hjelmslev planes are the main subject in \cite{LorimerF1,LorimerF2,LoLa2} while \cite{Lo1,Lo1b} deal with the structure of Hjelmslev rings.

\section{The continuation of the Hjelmslev epoch under Drake, Jungnickel and Sane} 

In his publications on $n$--uniform planes, we can observe already that Drake had a particular interest in {\sl finite} Hjelmslev planes. This is continued when his attention goes more and more to the problem of existence and non--existence of finite Hjelmslev planes with given parameters. In a series of papers \cite{Drake9,Drake10,Drake11,Drake12,DrakeJung3,DrakeL2,DrakeT}, some of them with co--author, he attacked this problem  and he linked finite PH--planes to nets in \cite{DrakeH,DrakeSH}.\\ Meanwhile, also finite Klingenberg planes came to the attention \cite{Drake13}. Drake and Lenz considered a parameter set for finite PK--planes in \cite{DrakeL1} together with new examples of finite PH--planes.\\ Structure theorems for finite chain rings (needed for finite desarguesian Klingenberg planes) were proved by Edwin {\sc Clark} {\sl et al.} in \cite{ClarkLiang, ClarkDrake} and independently by Arnold {\sc Neumaier} in \cite{Neu} and Al--Khamees \cite{AlKhamees}. The classification of all chain rings is still an open problem but partial results are known. Galois rings GR($q^n,p^n$) with $q^n$ elements and characteristic $p^n$, with $q=p^r$, play a crucial role. 

Beside Drake another player came to the forefront, Dieter {\sc Jungnickel}, who was active at the Universities of Giessen and Augsburg (Germany).
He was a student of Hanfried Lenz and became an expert in the theory of designs. In 1976 he wrote his Diplomarbeit at the University of Berlin (Germany) on {\sl ``Klingenberg and Hjelmslev planes"} and with the dissertation {\sl ``Konstruktion transitiver Inzidenzstrukturen mit Differenzenverfahren"} he obtained his doctoral degree. His most important contribution to the theory of Hjelmslev planes (and the more general class of Klingenberg planes) concerns \mbox{``regularity"} \cite{DrakeJung4,Jung2,Jung3,Jung5,Jung8}. A PK-- or PH--plane is regular if it has an abelian automorphism group $G = Z \oplus N$, where $G$ acts regularly (sharply transitively) on the point set and on the line set and $N$ acts regularly on each neighborhood. It is proved in \cite {HaleJung} that any finite PK--plane over a commutative local ring is regular. Regularity is also interpreted in terms of difference sets and auxiliary matrices, leading to new families of finite Hjelmslev and Klingenberg planes. An interesting result, connecting PH--planes with designs, is: the projective uniform Hjelmslev planes of order $q$ (with $q>2$) are precisely the symmetric divisible partial designs on two classes with parameters $v=b=q^2(q^2+q+1),k=r=q(q+1),s=q^2,t=q^2+q+1,\lambda_1 =q, \lambda_2 = 1$. For $q=2$ counterexamples exist, see \cite{Jung6}. The concepts of regularity and uniformity were also considered in $K$--structures, a further generalization of Klingenberg planes (see \cite{Drake13,DrakeJung1,DrakeJung2,Jung1,Jung4,Jung7,Jung9,Jung10}). Nino {\sc Civolani} \cite{Civolani} considers free extensions of partial Klingenberg planes. 

Jungnickel's work was continued by Sharad {\sc Sane}. Sane studied at the Indian Institute of Technology Bombay, Mumbai (India) and obtained his PhD with the dissertation {\sl ``Studies in Partial Designs and Projective Hjelmslev Planes"} under the supervision of Balmohan Vishnu {\sc Limaye}. In \cite{LimSane} Sane and Limaye demonstrate that $n$--uniform PH--planes are a kind of divi\-sible partial designs and by taking advantage of this property, they can give an alternative proof for the fact that $n$--uniform planes are strongly $n$--uniform, as was proved before in another way by T\"orner \cite{Torner2}.
In some other papers \cite{DrakeS1,DrakeS2,Sane1,Sane2,Sane3,Sane4,Sane5} Sane contributes to the theory of finite Hjelmslev and Klingenberg planes.

\section{The coordinatization of Hjelmslev and Klingenberg planes: a versatile story}

The coordinatization of affine and projective planes is one of the most power\-ful tools in the study of such geometries. It permits to reformulate  geometric properties into algebraic ones (and vice versa), leading to a better insight, including the construction of many non--desarguesian examples. This coordinatization goes back to Marshall {\sc Hall} Jr. (1910--1990). His important paper \cite{Hall} published in 1943 is still one of the most cited. The basic concept is a Hall ternary ring, also called PTR (planar ternary ring), an algebraic structure $(R,T)$ with $R$ a non--empty set containing two distinct elements $0$ and $1$ and with $T$ a ternary operation on $R$ such that $y=T(x,m,k)$ means that the point with coordinates $(x,y)$ lies on the line with coordinates $[m,k]$. With $(R,T)$ one can associate two loops $(R,+)$ and $(R,\circ)$ when \mbox{defining} $a+b:=T(a,1,b)$ and $a \circ b:=T(a,b,0)$. The properties of the plane (formulated in terms of the validity of Desargues' configuration or in terms of transitivity of the automorphism group) are reflected in the richness of the coordinatizing algebraic structure. If the theorem of Desargues is always valid or equivalently if the plane is $(P,\ell)$--transitive for any choice of the point $P$ and the line $\ell$, then it turns out that $(R,+)$ and $(R\setminus\{0\},\circ)$ are both groups, hence $(R,+,\circ)$ is a division ring or skewfield. Conversely, any skewfield gives rise to a desarguesian projective plane.
Slight variations on Hall's coordinatization method were made by Daniel {\sc Hughes} \cite{Hughes,HughesPiper} and by G\"unter {\sc Pickert} \cite{Pickert}. Independently the russian mathematician Lev Anatolevich {\sc Skornyakov} described in 1949 a similar coordinatization method in \cite{Skor1}. His work \cite{Skor2} was important for the distribution of the knowledge on projective planes in the russian speaking mathematical community.

In the seventies and the eighties several attempts were made to coordinatize in a similar way affine and projective Hjelmslev and Klingenberg planes. 
In 1967 the russian geometer V.K. {\sc Cyganova} worked out a first successful coordinatization for the more restrictive class of affine Hjelmslev planes \cite{Cyg2}. She used the concept of an {\sl $H$--ternar}, an algebraic structure with two ternary opera\-tions, generalizing a Hall ternary ring (one of the main differences being the existence of zero divisors). 
Since her paper was written in russian, it remained unfortunately unaccessible for many people. \\
Independently from Cyganova,  J.W. {\sc Lorimer} introduced in 1971 in his PhD thesis ``{\sl Hjelmslev Planes and Topological Hjelmslev Planes}" genera\-lized ternary rings (very similar to $H$--ternars) as the coordinatizing structures of affine Hjelmslev planes. Three years later, Phyrne {\sc Bacon} streamlined the work of Cyganova and Lorimer in her thesis ``{\sl Coordinatized Hjelmslev planes}", and she introduced the name {\sl biternary ring} (in appendix A of her thesis she gives a comprehensive list of annotations including some mistakes and imperfections in the work of her predecessors). The interaction between the geometric properties of an affine Hjelmslev plane and the algebraic proper\-ties of its coordinatizing biternary ring was examined in more detail by Lorimer in \cite{LorimerC1}, by Cyganova in \cite{Cyg1,Cyg3,Cyg4,Cyg5,Cyg6,Cyg7,Cyg8}, by {\sc Emelchenkov} in \cite{Emelchenkov1,Emelchenkov3,Emelchenkov4} and by {\sc Shatokhin} in \cite{Shatokhin1,Shatokhin2}.\\ To be complete, we also have to mention a paper by Drake \cite{DrakeC} in which he obtains a kind of coordinatization for a special class of affine Hjelmslev planes (radial planes) by means of a module.

After the coordinatization of affine Hjelmslev planes, a similar theory for the more general class of affine Klingenberg planes has been worked out by several authors. {\sc Bacon} generalized her biternary rings. Drake, her supervisor, encouraged her to publish in mathematical journals, but she was stubborn and, apart from one single publication \cite{Bacon4}, she refused. Her voluminous work, totaling about 1000 pages, is contained in four books \cite{Bacon3}, edited in own management. For that reason it was often overlooked and seldom recognized as an acceptable reference. In \cite{Bacon4} the ``triangle theorem" is proved: a PK--plane $\mathcal{P}$ possessing a nondegenerate triangle with sides $\ell_1, \ell_2$ and $\ell_3$ such that each derived AK--plane $\mathcal{A}_i=\mathcal{P}\setminus g_i$ is desarguesian, is itself a desarguesian PK--plane.\\
The Czech mathematician Frantisek {\sc Machala} from the University of Olomouc introduced in \cite{MachalaC3} {\sl affine local ternary rings} as an alternative for the coordinatization of affine Klingenberg planes. Much later he could prove the equivalence between his coordinatization method and the one of Bacon. He also proved that any ``incomplete"\ biternary  ring (with one ternary and one partial ternary operator) can be extended to a biternary ring with two (full) ternary operators \cite{MachalaC6}. In \cite{BakerLorimer} it is shown that this biternary ring extension is unique. 

In the case of ordinary planes each (desarguesian) affine plane can be extended to a (desarguesian) projective plane. This does not hold any longer for Hjelmslev planes (see section 7). This observation has also a serious impact on the coordinatization of projective Hjelmslev planes: it does not follow immediately from the affine coordinatization. 
The projective case was first attacked by the russian mathematician E.P. {\sc Emelchenkov} in 1972 in his PhD thesis {\sl "Ternars and automorphisms of Hjelmslev planes"} (in russian) and in \cite{Emelchenkov2}. Due to the language barrier, his work was not accessible for many researchers and for that reason, like Cyganova's work, it was somewhat overlooked.\\ 
Coordinatization methods for the more general case of projective Klingenberg planes were worked out by both {\sc Machala} and {\sc Bacon}. Their approach is totally different. Machala's method is based on the concept of an {\sl extended local ternary ring}, an algebraic structure $(R,R',T)$ with two disjoint sets of coordinates $R$ and $R'$ and one ternary operation $T$ on $R \cup R'$ (see \cite{MachalaC1,MachalaC2,MachalaC4}). This coordinatization was not very successful because it is not obvious to see any interaction between properties of the extended local ternary ring and geometric properties of the coordinatized plane.\\
The approach of Bacon is based on the fact that a PK--plane can be \mbox{covered} by three AK--planes corresponding to three biternary rings. This yields a coordinatizing structure for a projective Klingenberg plane as a triplet of biternary rings, called {\sl sexternary ring} \cite{Bacon3}. In Bacon's voluminous work, the interaction between geometric properties and the algebraic structure is examined in depth. Unfortunately, a big part of these results remained hidden for the reason mentioned above. In \cite{DugasC3} Manfred {\sc Dugas} used a similar coordinatizing structure, with six ternary operations.\\
Independently from the people mentioned above, the author introduced in 1987 in his PhD thesis {\sl ``Klingenberg incidence structures, a contribution to ring geometry"} (in Dutch) (see also \cite{KeppensC1,KeppensC2,KeppensC3}) {\sl planar sexternary rings (PSR's)} with one full and five partial ternary operators. His coordinatization method for PK--planes was inspired by the Hughes variant of the Hall ternary ring. A small deficiency in his method was detected later (as pointed out by Baker and Lorimer in {\cite{BakerLorimer}).  
As a consequence of this shortcoming, the coordinatization of a PK--plane by a PSR wasn't fully compatible with the coordinatization of a derived AK--plane by the biternary ring obtained from the PSR. 
To overcome this anomaly, Baker and Lorimer (op.cit.) \mbox{developed} a new coordinate ring, called {\sl incomplete sexternary ring} (in the spirit of Dugas) as a substitute for the planar sexternary ring. They even proved that such a structure can be extended (in a unique way) to a {\sl sexternary ring} with six full ternary operators.

The coordinatization of projective planes is a handy instrument for the construction of non--desarguesian examples in an algebraic manner. A lot of new planes were found using quasifields, nearfields or alternative division rings. Because of the bigger complexity of sexternary rings it seems that much less examples of non--desarguesian PK--planes were obtained in this manner. Never\-theless examples of non--desarguesian AK-- and PK--planes obtained from algebraic structures, are known. The oldest examples are the Moulton affine Hjelmslev planes, given by Baker in \cite{Baker1}. A projective version is contructed in \cite{KeppensC3} by the author. Klingenberg planes over local nearrings and Hjelmslev planes over Hjelmslev--nearrings were studied by Emanuel {\sc Kolb} in \cite{Kolb1,Kolb2}. Much attention has gone to Moufang planes which can be coordinatized by local alternative rings. Moufang--Hjelmslev planes first appear in a paper of Dugas \cite{DugasC1}.   
He proves that all finite uniform Moufang PH--planes are desarguesian. A stronger version of that theorem was later proved in \cite{DugasC2} (the uniformity condition could be dropped if the order of the plane is bigger than 2). Similar results were found by Baker, Lane and Lorimer for Moufang Klingenberg planes, see \cite{BLL1,BLL3,BLL4}. They prove that the class of \mbox{Moufang} PK--planes coincides with the class of planes over local alternative rings and that a finite Moufang PK--plane in which any two points have at least one joining line, is a desarguesian projective Hjelmslev plane. Also a stronger version of Bacon's triangle theorem was proved in \cite{BLL2}: a PK--plane with a non--degenerate triangle for which the three derived AK--planes are translation AK--planes (and with epimorphic image distinct from PG(2,2)), is Moufang. \\More recently, a group of mathematicians around Basri {\sc \c{C}elik} and S\"uleyman {\sc \c{C}ift\c{c}i}, from the University of Uludag, Bursa (Turkey), published \mbox{several} papers concerning a particular class of Moufang--Klingen\-berg planes \cite{ACC1,ACC2,ACC3,ADDBG,CelikCiftci,ACC4,ACC5,ACC6,ACC7}. Their results, all variations on the same theme, overlap with work of Andrea {\sc Blunck} \cite{BlunckC,BlunckC2}.
The role of Pappus' theorem (its validity in a desarguesian plane implies the commutativity of the coordinatizing ring) was investigated by {\sc Nolte} and {\sc Maurer} in \cite{MaurerNolte,Nolte}.
 
\section{Order and topology in Hjelmslev geometry:  Machala and Lorimer}

The theory of ordered incidence structures can be traced back mainly to {\sc Pasch}, {\sl Vorlesungen \"uber neuere Geometrie} from 1882. An excellent survey paper on the axiomatics of ordered incidence geometry is \cite{Pambuccian}.\\
Ordered affine and projective Hjelmslev planes were studied by a group of Canadian mathematicians, starting in the seventies. At least three dissertations were written in that period at the McMaster University of Hamilton, Ontario (Canada) under supervision of Norman {\sc Lane} who was also a world--class canoeist, competing in two Olympic games (bronze medal in 1948 in London), before he started his academic career. \\In 1975, Lynda Ann {\sc Thomas}, wrote her Master's thesis on {\sl ``Ordered desar\-guesian affine Hjelmslev planes"} in which she proved that any ordered AH--ring gives rise to an ordered desarguesian affine Hjelmslev plane and vice versa. This result was published a few years later in \cite{Thomas}.\\ James {\sc Laxton}, another student of Lane, treated in his Master's thesis {\sl ``Ordered non--desarguesian affine Hjelmslev planes"}.
Catherine {\sc Baker}, also a student of Lane, wrote her Master's thesis on {\sl ``Affine Hjelmslev and generalized affine Hjelmslev planes"} and her doctoral thesis in 1978 on {\sl ``Ordered Hjelmslev planes"}. In that thesis she investigates in detail the relationship between ordered AH--planes and the coordinatizing ordered biternary rings,  extending results of Laxton and Thomas. Baker \mbox{published} several papers about ordered (affine and projective) Hjelmslev planes (also with co--author): \cite{Baker2,Baker3,Baker4,Baker5,Baker6,Baker7}. \\ It would be disrespectful if we should attribute all the results on ordered ring geometries to the ``Canadian School" only. Independently,  a theory of orderings for Klingenberg planes was worked out by {\sc Machala}. He published many papers on ordered Klingenberg planes \cite{Machala1,Machala2,Machala3,Machala4,Machala5,Machala6,Machala7,Machala8} and one overview work \cite{Machala9}. The work of Baker et al.~resembles in many aspects Machala's, but there are some differences. For a comparison between both approaches, one may consult \cite{Baker7}.

The study of {\sl topological} Klingenberg and Hjelmslev planes remained an exclusive Canadian affair. 
The most prominent student of Lane, was undoubtedly J.W. (Mike) {\sc Lorimer}. In his doctoral thesis on {\sl ``Hjelmslev planes and topological Hjelmslev planes"}, he not only introduced a coordinatization (see the previous section) but he also laid the foundation for topo\-logical Hjelmslev geometry. His work generalizes that of Salzmann \cite{Salzmann} and \mbox{Skornyakov} \cite{Skornyakov} on topological projective planes. In a series of publications \cite{Lorimer1,Lorimer2,Lorimer3,Lorimer4,Lorimer5,Lorimer6,Lorimer7,Lorimer8,Lorimer9,Lorimer10,Lorimer11} he further developed the theory in close connection with the coordinatization problem. Among the most important theorems proved by Lorimer, we mention the following characterization theorem: the only locally compact connected pappian projective Hjelmslev planes are the ones over the rings $\mathbb{K}[x]/\langle x^n \rangle$ with $\mathbb{K}$ the field of real or complex numbers.

\section{The revival of ring geometry in the eighties and nineties under Veldkamp and Faulkner}

In the seventies ring geometry was restricted almost exclusively to Hjelmslev and Klingenberg geometry (in the desarguesian case to geometries over local rings and Hjelmslev rings).
The Dutch mathematician Ferdinand Douwe {\sc Veldkamp} (1931-1999), who is well--known for his work on geometries associated with exceptional Lie groups and in particular polar spaces, reverted back to the pioneering work of Barbilian where geometries over the broader class of $Z$--rings were considered. It was Veldkamp's aim to give an axiom system for projective planes (and higherdimensional spaces) over arbitrary rings with unit (without the imperfection in Barbilian's attempt \cite{Barbi1,Barbi2}). From conversations with his colleague van der Kallen at the University of Utrecht (an expert in $K$--theory), it became clear that the best setting for this project is provided by rings of stable rank two. A ring $R$ has stable rank two if the following property holds: if $a,b \in R$ and $Ra+Rb=R$ then there exists a $r$ in $R$ such that $a+rb$ is invertible in $R$. The class of stable rank two rings comprises all local rings. Hence, the projective ring planes introduced by Veldkamp include the desarguesian Klingenberg and Hjelmslev planes. A ring of stable rank two is always a $Z$--ring. \\Veldkamp first worked out the theory for planes in \cite{Veldkamp1} with some special cases in \cite{Veldkamp2} and later for spaces of higher dimension (see section 12). 
The ring planes defined by Veldkamp are also known today as desarguesian Veldkamp planes. 
 
John Robert {\sc Faulkner}, an authority in the domain of non--associative algebra and geometry, further extended the theory of Veldkamp planes in the non--desarguesian direction by introducing alternative (non--associative) rings of stable rank two in \cite{Faulkner1}. He then proved in \cite{Faulkner2} that a Veldkamp plane has the Moufang property (i.e.~$(P,\ell)$--transitivity holds for all $P,\ell$ with $P$ incident with $\ell$) if and only if it is a plane $\mathcal{P}(\mathbb{A})$ over an alternative ring $\mathbb{A}$ of stable rank two. Inspired by the work of his predecessors, Faulkner introduced in \cite{Faulkner3} Faulkner planes as a very general class of plane incidence structures with neighbor (or remoteness) relation. They comprise the Veldkamp planes and the planes introduced by Barbilian. A (connected) Faulkner plane for which the group of $(P,\ell)$--transvections (automorphisms fixing all objects incident with $P$ or $\ell$) is transitive on the set of points not neighboring with $\ell$, is called a transvection plane. The coordinatization of transvection Faulkner planes by a not necessarily associative alternative ring with the property that $ab=1$ implies $ba=1$ involves a rather technical procedure based on group theory and a lot of new concepts such as covering planes and tangent bundle planes. A transvection Faulkner plane for which the tangent bundle plane is also a transvection plane is called a Lie transvection Faulkner plane. To every such plane an alternative two--sided units ring can be attached and conversely with an alternative two--sided units ring $R$ a corresponding Lie transvection Faulkner plane can be constructed. However, this plane is not determined unambiguously when the ring does not have stable rank 2. This high price has to be paid for the generalization from Veldkamp planes to Faulkner planes. In \cite{Faulkner4} Faulkner gives a geometric construction of Barbilian planes coordinatized by composition algebras (including the Moufang plane) using  Jordan algebras. His book \cite{Faulkner7} is completely devoted to the role of such algebras in projective geometry.\\
Faulkner was surrounded by students at the University of Virginia, Charlottesville (USA), who all were involved with the study of ring geometries. Terese Deltz {\sc Magnus} obtained her PhD in 1991 on {\sl ``Geometries over non--division rings"} and  she could generalize Faulkner's axioms and results for geometries of higher dimension (see section 12). Eve {\sc Torrence} also graduated under Faulkner's supervision with {\sl ``The coordinatization of a hexagonal--Barbilian plane by a quadratic Jordan algebra"}, a generalization of the classical notion of generalized hexagon. Another student was \mbox{Catherine} {\sc Moore d'Ortona} who studied homomorphisms between projective ring planes in her PhD thesis {\sl ``Homomorphisms of remotely projective planes"}, published in \cite{Moore}. Finally Karen {\sc Klintworth} wrote her PhD thesis on {\sl ``Affine remoteness planes"}.

Faulkner himself considered a slightly more general axiomatization of Faulkner planes in \cite{Faulkner6}. There he chooses for the remoteness relation, the negation of the neighbor relation. Much of the results obtained in \cite{Faulkner3} are extended but the coordinate rings that appear are no longer always alternative. It is proved that $\mathcal{P}$ is a transvection plane if and only if $\mathcal{P}$ is isomorphic to $P(G,N)$, the plane associated with a group $G$ of Steinberg type parametrized by the ring $R$ and with $N$ a certain subgroup of $G$. Necessary and sufficient conditions are given for $R$ to be alternative, associative or commutative. In \cite{Faulkner6} also projective remoteness planes with reflections (hence metric planes) have been considered.  
The content of this paper is closely related to some work of Kn\"uppel and Salow in \cite{Knuppel5} (see section 6). It also contains a part on affine ring planes and elementary basis sets which are closely related to Barbilian domains as introduced by Leissner (see section 5). 

In the slipstream of Veldkamp's paper on projective ring planes, several slightly modified axiom systems have been described, leading to other classes of projective ring geometries.
We have already mentioned Frieder {\sc Kn\"uppel} in the section on metric ring geometry, but some of his papers rather join the spirit of this section. In \cite{KnuppelE1} Kn\"uppel considers ring geometries over associative rings based on four axioms adapted from Veldkamp (based on remoteness rather than neighborship). A coordinatization theorem is stated without proof. In \cite{KnuppelE2} he studies homomorphisms between such planes.\\
Renata {\sc Spanicciati} defines near--Barbilian planes and strong near--Barbilian planes in \cite{Span} by adapting some of the axioms of Veldkamp. The neighbour relation between points turns out to be the identity, and the neighbour relation between lines becomes an equivalence relation. In \cite{HansVM} Guy {\sc Hanssens} and Hendrik {\sc Van Maldeghem} prove that any near--Barbilian plane is strong near--Barbilian. K\'alman {\sc P\'entek} gives a necessary and sufficient condition for a Veldkamp plane to be a direct product of a finite number of Veldkamp planes in \cite{Pentek1,Pentek2,Pentek3}.\\
The study of homomorphisms between ring geometries was also a central theme in several papers of Veldkamp (some of them in joint work with Joseph {\sc Ferrar}) \cite{Ferrar,Veldkamp6,Veldkamp5,Veldkamp3,Veldkamp4}.
The geometric homomorphisms of distinct kind (incidence preserving, neighbor--preserving, distant--preserving) are characterized in terms of algebraic morphisms between the underlying rings.
Veldkamp's results were generalized by the author in the non--desarguesian case for homomorphisms between projective Klingenberg planes using the coordinatizing planar sexternary rings (see section 9) in \cite{KeppensH1,KeppensH2}. \\Thorsten {\sc Pfeiffer} could generalize a well--known theorem for planes over fields to planes over rings: a desarguesian ring plane $\mathcal{P}(R)$ is pappian (Pappus' theorem is valid) if and only if $R$ is commutative. 

\section{Projective and affine Hjelmslev spaces and spaces over arbitrary rings}

Hitherto we only discussed {\sl plane} ring geometries. The theory of higher dimensional projective spaces over rings has been developed by different people.
Projective spaces over local rings appear for the first time in the work of {\sc Klingenberg} \cite{Kling4}. Today we call them Klingenberg spaces (PK--spaces).
The first study of PK--spaces after Klingenberg, is due to Hans--Heinrich {\sc L\"uck}, a student of L\"uneburg. His PhD thesis, published as an article \cite{Luck} in 1970 under the somewhat misleading title {\sl ``Projektive Hjelmslevr\"aume"},  contains an axiomatic characterization of a class of incidence structures which permit a coordinatization by local rings. Hence, the paper deals with projective Klingenberg spaces rather than with Hjelmslev spaces.
L\"uck proves that in an axiomatically defined PK--space of dimension at least three, the theorem of Desargues holds and that it must be isomorphic to a space derived from a module over a local ring. Hence all projective Klingenberg spaces of dimension $\geq 3$ are desarguesian, a situation analogous to the case of classical projective spaces. \\
Independently from L\"uck, {\sc Machala} defined and \mbox{studied} projective Klingenberg spaces (Projektive R\"aume mit Homomorphismus) of finite or infinite dimension in \cite{Machala10, Machala11, Machala14, Machala12, Machala13}. He proved that the planes in a PK--space are PK--planes and that PK--spaces of dimension at least three come from modules over a local ring (cfr. L\"uck). Machala also investigated homomorphisms between PK--spaces and the fundamental theorem (isomorphisms between spaces can be represented by semilinear mappings between the underlying modules). A paper by Jukl \cite{Jukl2} is in conformity with this.\\
An axiomatic approach for the more restricted class of projective Hjelmslev spaces (PH--spaces) was initiated by John {\sc Lamb} Jr. from the University of Texas at Austin (USA) in his PhD thesis, entitled {\sl ``The Structure of Hjelmslev space, a generalization of projective space"} (1969) but its content (related to lattice theory) was not published. \\
Much more widespread is the work of Karl {\sc Mathiak}, who was very productive in this field. He defined a class of special projective Hjelmslev spaces \mbox{starting} from a vectorspace over a (skew)field endowed with a valuation (``Bewertete Vektorr\"aume"). The structure of the ideals in the corresponding \mbox{valuation} ring plays a central role in his approach. The theory is \mbox{thoroughly} worked out in a series of papers, published over a period of twenty years between 1967 and 1987, see \cite{Mathiak1,Mathiak2,Mathiak3,Mathiak4,Mathiak5,Mathiak6,Mathiak7,Mathiak8,Mathiak9,Mathiak10}.\\ His compatriot Alexander {\sc Kreuzer} introduced an \mbox{axiom} system for arbitrary PH--spaces in his doctoral thesis {\sl ``Projektive Hjelmslevr\"aume"} which was published afterwards as an article in \cite{Kreuzer2} with some preliminary work in \cite{Kreuzer1}. This study is continued in \cite{Kreuzer3,Kreuzer4}. 

For the affine case we have to go to Canada again. In four papers, Tibor {\sc Bisztriczky} together with J.W. {\sc Lorimer}, worked out two axiom systems for affine Klingenberg spaces \cite{Bis1,Bis2,Bis3,Bis4}. Neither of their axiom systems assumes the existence of an overlying projective Klingenberg space or the existence of an underlying ordinary affine space. Machala in \cite{Machala12} also defined affine Klingenberg spaces, but not separate from PK--spaces (similar to ordinary affine spaces obtained from projective spaces by deleting a hyperplane). 

The most general study of projective ring spaces (over not necessarily local rings) was undertaken by Ferdinand {\sc Veldkamp}. We have already indicated his interest in section 11. In \cite{Veld1,Veld2}, Veldkamp gives a self--dual axiom system for projective ``Barbilian spaces" of finite dimension using the basic concepts of points, hyperplanes, incidence and neighbor relation. We call them now Veldkamp spaces. The main result is that Veldkamp spaces of dimension $\geq 3$ are spaces over rings of stable rank 2. \\ Theresa {\sc Magnus} defines Faulkner spaces as even more general geometries by extending the theory of Faulkner planes \cite{Magnus}. She proves that any Faulkner space of dimension $n \geq 3$ is coordinatized by a unique associative two--sided units ring $R$ and that the group generated by all transvections is a group of Steinberg type over $R$. A Faulkner space over the ring $\mathbb{Z}$ of integers is constructed, providing an example of a Faulkner geometry which is not a Veldkamp space, since $\mathbb{Z}$ has stable rank 3. Under the Veldkamp spaces we also find the projective spaces over matrix rings over GF($q$), studied by {\sc Thas} \cite{Thas3} in the early days of ring geometry (see section 4). Other contributions to finite ring spaces are due to {\sc Kapralova} who considered projective spaces over the ring of dual numbers over a Galois field in \cite{Kapralova} and to Ivan {\sc Landjev} and Peter {\sc Vandendriesche} \cite{LVD1,LVD2}.

A well-developed theory for affine spaces over rings is also due to {\sc Veldkamp}. In \cite{Veld3} he defines Barbilian domains in free modules of rank $n$ and introduces $n$--dimensional affine ring geometries. A geometrical interpretation of Barbilian domains is given by Sprenger in \cite{Sprenger}.\\
Other attempts for setting up a theory of higher dimensional affine ring geometries (incidence structures with parallellism) are scattered in the literature. The definitions and the methods used are very diverse. Contributions in this field are due {\sc Permutti} and {\sc Pizzarello} \cite{Permutti,Pizzarello}, {\sc Miron} \cite{Miron}, {\sc Leissner, Severin} and {\sc Wolf}  \cite{Leis4,LSW}, {\sc Ostrowski} \cite{Ostrowski}, {\sc Schmidt} and {\sc Weller} \cite{SW}, {\sc Kreis} \cite{Kreis,KreisSchmidt}, {\sc Seier} \cite{SeierS1,SeierS2}, {\sc Bach} \cite{Bach} and others.

One of the problems for higher dimensional geometries over rings that has got much attention is that of morphisms and the fundamental theorem. For classical projective geometries $P(V)$ induced by a vectorspace over a field this theorem states that any bijective incidence preserving map (projectivity) between projective spaces $P(V)$ and $P(W)$ can be algebraically characterized by a semi\-linear map from $V$ to $W$.
The first generalization of the fundamental theorem to ring geometries was obtained by the Indian mathematicians {\sc Ojanguren} and {\sc Sridharan} who could prove it in case of module--induced geometries $P(M)$ with $M$ a free module of finite rank $\geq 3$ over a commutative ring \cite{OS}. Generalizations to other classes of rings were proved later by {\sc Sarath} and {\sc Varadarajan} in \cite{SV} and by {\sc James} in \cite{James} and {\sc Faure} in \cite{Faure}. \\
For the sake of completeness we also refer to module--induced geometries as defined by Marcus {\sc Greferath} and Stefan {\sc Schmidt} \cite{Gref1,Gref2,Gref3,Gref4}. Instead of the usual definition of the pointset as the set of all submodules generated by a unimodular element in a free module (cfr. Veldkamp), they take all submodules of rank one, leading to the bizarre situation of points properly contained in bigger points. In \cite{Gref0} an extension of the fundamental theorem is proved for such module--induced geometries. Close relative to this, there is an abundance of articles by the school of the Georgian mathematician Alexander {\sc Lashkhi}. They all contain variations on the same theme:  an extension of the fundamental theorem for affine and projective geometries related with modules over rings, from the lattice--theoretic point of view. We have not included the whole collection of papers by Lashki and his students. Some of them have been published multiple times in different journals (in Russian and in English). The literature list only mentions a few representative ones:  \cite{Kviri1,Kviri2,Kviri3,Lash1,Lash2,Lash3,Lash4,Lash5,Rost}.  

\section{Projective lines and circle geometries over rings: Blunck, Havlicek and Keppens}

The ``smallest" projective geometries that can be considered, are the {\sl projective lines}. From the viewpoint of incidence geometry not much can be said about these rather poor structures. But combining the study of projective lines with those of their automorphisms (the general linear group) puts them into a new light. The theory has common ground with what is  usually called {\sl geometric algebra}. Indeed, the projective line $P(R)$ over any ring $R$ can be defined in terms of the free left $R$--module $R^2$ as follows: it is the orbit of a starter point $R(1,0)$ under the action of the general linear group $GL_2(R)$ on $R^2$. Since geometric algebra over rings only fits sideways in the section of incidence geometry, we did not make an effort to be complete in the literature list for this item. Nevertheless we mention a number of relevant references in which more information can be found, e.g.~\cite{McDonald2}. Central themes that keep returning are the notions of cross--ratio and harmonic quadruples and the fundamental theorem, also known as Von Staudt's theorem. In the case of classical projective lines over a field or a skewfield, this theorem characterizes mappings of the projective line which preserve harmonicity as projectivities. 

Among the first publications on projective ring lines (and we do not consider here papers only dealing with linear groups over rings) belong some articles by the Indian mathematicians Nirmala and Balmohan {\sc Limaye}. They prove a generalization of Von Staudt's theorem for some special classes of commutative and non--commutative rings in \cite{Lim1,Lim2, Lim3,Lim4,Lim5}. A little bit earlier, in 1968, {\sc Melchior}, a student of Benz, wrote his PhD thesis, entitled {\sl ``Die projektive Gerade \"uber einem lokalen Ring: ihre lineare Gruppe und ihre Geometrie"}. Other contributions to this item appeared in \cite{BM2,BaBa,Bart,Blunck5,Blunck9,Chkhat,Cir,Havl1,Havl2,Herzer,Jurga,Kulkarni,KL,LashR5,McDonald1,Sanigaetal,Schaeffer}. Some papers by {\sc Bilo} and {\sc Depunt} \cite{Bilo}, {\sc Hubaut} \cite{Hubaut1,Hubaut2} and {\sc Thas} \cite{Thas1,Thas2} also deal with projective lines over rings (see section 4). They were followed by {\sc Havlicek} {\sl et al.} in \cite{Havl4,Havl5,Havl6}.

Projective lines over rings are also intrinsically related with circle geometries. This relation was established for the first time by {\sc Benz} in his famous book {\sl ``Vorlesungen \"uber Geometrie der Algebren"} \cite{Benzboek} from 1973. He presented a unified treatment of plane circle geometries, now called {\sl Benz planes}, using the projective line over a commutative ring which is a two--dimensional $\mathbb{K}$--algebra over a field $\mathbb{K}$. His definition was extended by Andrea {\sc Blunck}, a student of Benz, and Armin {\sc Herzer} who introduced more general chain geometries $\Sigma(\mathbb{K},R)$ with $R$ a (not necessarily two--dimensional) $\mathbb{K}$--algebra. A chain geometry is an incidence structure whose point set is the set of points of the projective line over $R$ and the $GL_2(R)$ orbit of $P(\mathbb{K})$ is the set of chains. The plane circle geometries of M\"obius, Laguerre and Minkowski type are particular chain geometries for $R=\mathbb{L}$ (a quadratic field extension of $\mathbb{K}$), $R=\mathbb{K}+\mathbb{K}\varepsilon$ with $\varepsilon^2=0$ (dual numbers) or $R=\mathbb{K}+\mathbb{K}t$ with $t^2=t$ (double numbers) respectively. For an overview of chain geometry we refer to \cite{BlunckHerzer}. Chain geometries were also treated by {\sc Schaeffer}, a student of Benz, in his doctoral thesis entitled {\sl ``Zum Automorphismenproblem in affinen Geometrien und Kettengeometrien \"uber Ringen}. The study of chain geometries is continued by {\sc Blunck} who introduced generalized chain geometries by considering projective lines over non--associative, alternative rings. In her PhD thesis {\sl ``Doppelverh\"altnisse und lokale Alternativringe"} (1990) and in \cite{Blunck1} she extended the notion of cross--ratio and investigated chain geometries (in relation to projective lines over non--associative rings) in \cite{Blunck3,Blunck6,Blunck6bis}.

Around the same time the author  in \cite{Keppens1,Keppens2} defined Klingenberg--Benz planes axiomatically, i.e.~plane circle geometries with neighbor relation which admit a natural epimorphism onto a classical Benz plane. Using the projective line over three kinds of quadratic ring extensions of a local ring (instead of a field), he was able to construct algebraic models of such geometries. Also Konrad {\sc Lang} in \cite{Lang} studied independently of us a class of Hjelmslev--M\"obius planes. {\sc Blunck} and {\sc Stroppel} extended our definition of Klingenberg--Benz planes to Klingenberg chain spaces in \cite{BlunckStroppel} and Blunck also proved that a Klingenberg chain space can be embedded into a projective Klingenberg space, such that the points are identified with points of a quadric and the chains with plane sections \cite{Blunck4}. In \cite {SeierC1} {\sc Seier}  constructed analogously a class of chain geometries $\Sigma(H,L)$ with $H$ a Hjelmslev--ring and $L$ a ring extension of $H$ and in \cite{SeierC2} he defined a M\"obius plane with neighbor relation of a different kind than the one defined by us.\\
A basic notion concerning the projective line over a ring $R$ is its distant relation: two points are called distant if they can be represented by the ele\-ments of a two-element basis of $R^2$. The distant graph has as vertices the points of the projective line and as edges the pairs of distant points. The distant graph is connected precisely when $GL_2(R)$ is generated by the elementary linear group $E_2(R)$. This aspect of projective ring lines was studied in more detail by {\sc Blunck}, {\sc Havlicek}, {\sc Matra\'s} and some others  in \cite{BM1,Blunck7,Blunck8,Blunck8bis,Havl3,MS,MS2}. In \cite{Blunck10,Blunck11} the interaction between ring geometry and the geometry of matrices in the sense of Hua (see \cite{Hua}) is investigated in more detail.  

\section{Ring geometries and buildings: Van Maldeghem and co.}

The theory of buildings was invented by the Belgium-born French mathematician Jacques {\sc Tits}. Roughly speaking, buildings are incidence geometries on which groups act. Tits received many awards for his fundamental and path--breaking mathematical ideas, including the Abel Prize in 2008.
One of his achievements is the classification of affine buildings of rank at least 4. They are known to be ``classical", i.e.~they arise from algebraic groups over a local field. In the rank three case (where affine buildings are of three possible types $\tilde{A}_2$, $\tilde{C}_2$ or $\tilde{G}_2$) many non--classical counterexamples are known. \\
In his PhD thesis {\sl ``Niet--klassieke driehoeksgebouwen"} on triangle buildings (affine buildings of type $\tilde{A}_2$) Hendrik {\sc Van Maldeghem} observed that a special kind of ring geometry is present as the suitably defined geometry at distance $n$ from any given vertex of the building (the so-called $n$--th floor). This was described among other things in \cite{VM1,VM2}. A little bit later {\sc Hanssens}  and {\sc Van Maldeghem} could prove that those ring geometries are in fact projective Hjelmslev planes of level $n$, see \cite{VM3}. In \cite{VM4} they give a universal construction for level $n$ Hjelmslev planes (see also \cite{HallRao} for the $2$--uniform case) and as a corollary any level $n$ projective Hjelmslev plane is isomorphic to the $n$--th floor of a triangle \mbox{building}. This result links the theory of PH--planes to that of triangle buildings, a rather unexpected but fascinating fact. In the same spirit Van Maldeghem investigated another class of rank three affine buildings, of type $\tilde{C}_2$, and proved that the $n$--th floor turns out to be another type of ring geometry which can be seen as a generalization of an ordinary generalized quadrangle. He called it ``Hjelmslev--quadrangle" of level $n$ (see \cite{VM7,VM8}). In joint work with Hanssens a complete characterization of $\tilde{C}_2$--buildings by Hjelmslev quadrangles was obtained \cite{VM5,VM6}.\\ The author defined ``Klingenberg--quadrangles" as another generalization of ordinary generalized quadrangles in \cite{Keppensx}. The connection between Klingen\-berg--quadrangles and Hjelmslev--quadrangles is explained in \cite{VM7}. 

The relation of polar spaces of higher rank to generalized quadrangles is comparable with the relation of projective spaces to projective planes. Generalized quadrangles are polar spaces of rank two.
Projective ring spaces of dimension at least three have got as much attention in the literature as projective ring planes. This is not the case so far for general ``Klingenberg--polar spaces" or ``polar spaces over rings" if the rank is bigger than two.  Only one paper by {\sc James} \cite{JamesP} is known to us. Certainly this topic offers perspectives for future research.

The discovery by Van Maldeghem of the connection between buildings and ring geometry has a precedent. Twenty years earlier, in 1968, {\sc Veldkamp} in joint work  with Tonny {\sc Springer} con\-sidered a geometry over the split octonions (over the complex number field) in \cite{VeldSpringer}. This geometry is a kind of analogue of the non--desarguesian projective plane over the alternative division ring of (non--split) octonions $\mathbb{O}$ but in which two distinct lines may have more than one point in common and dually. It was called Moufang--Hjelmslev plane, but this name is misleading since it is not a projective Hjelmslev plane (the neighbor relation is not transitive) and hence it is completely distinct from the Moufang PH--plane (over an alternative local ring) studied elsewhere (see section 9). In two subsequent papers \cite{VeldHM1,VeldHM2} more results on Hjelmslev--Moufang planes are obtained, \mbox{concerning} projective groups.  The geometry of Veldkamp and Springer is the same as the one constructed by Tits \cite{Tits} starting from the split algebraic group of type $E_6$. In this geometry each line has the structure of a polar space and two lines can meet in more than one point (namely, in a maximal singular subspace of the corresponding polar spaces).\\ John {\sc Faulkner} considered Hjelmslev--Moufang planes by allowing an arbitrary ground field instead of $\mathbb{C}$ in \cite{Faul1,Faul2} and also Robert {\sc Bix} studied generalized Moufang planes in \cite{Bix1,Bix2}.\\
The relationship between Hjelmslev planes and buildings was further exploited by {\sc Van Maldeghem} and {\sc Van Steen} to give a characterization of some rank three buildings by automorphism groups \cite{VS1,VS2,VS3}. 

In the margins of the study of buildings some other questions have emerged. One of them concerns embeddings. Embeddings of point--line geometries into projective spaces are well--known in the literature. The embedding question for ring geometries, in particular for projective Hjelmslev planes, is first attacked by {\sc Artmann} \cite{Art4} who shows that the PH--plane over the ring of plural numbers $\mathbb{F}[t]/t^n$ ($\mathbb{F}$ a field), can be embedded in the $(3n-1)$--dimensional projective space over $\mathbb{F}$. In \cite{KeppensVM} the author and {\sc Van Maldeghem} prove a nice characterization theorem for embeddable Klingenberg planes: if $\mathcal{P}$ is a projective Klingenberg plane that is fully embedded in the projective space PG(5,$\mathbb{K}$) for some skewfield $\mathbb{K}$, then $\mathcal{P}$ is either a desarguesian Klingenberg plane over a ring of twisted dual numbers or a subgeometry of an ordinary projective plane. The embedding of the projective plane over a matrix ring with entries in GF($q$) into a projective space over GF($q$) was also observed by {\sc Thas} in \cite{Thas3,Thas4}.\\
Veronesean sets are closely connected with embeddings. In \cite{SVM1} {\sc Schillewaert} and {\sc Van Maldeghem} define geometries with an additional axiom by which the Hjelmslev--Moufang plane (in the sense of Springer--Veldkamp) and its relatives fit into the framework using the modern notion of para\-polar spaces. In \cite{SVM2} they provide a common characterisation of projective planes over two-dimensional quadratic algebras (over an arbitrary field) in terms of associated Veronesean sets. Anneleen {\sc De Schepper} and {\sc Van Maldeghem} \cite{DSVM} have considered Veronese representations of Hjelmslev planes over quadratic alternative algebras as part of a more general study of Veronese varieties and Mazzocca--Melone sets. 

\section{Ring geometries in coding theory: Honold, Kiermaier and Landjev}

One of the fastest growing disciplines in mathematics is coding theory. Since its introduction by Claude {\sc Shannon} in 1948, the number of publications about codes has exploded, in particular due to its importance in cryptography, data transmission and data storage. Initially mostly linear codes over finite fields were studied, but after the publication in 1994 of the paper \cite{Hammons} by {\sc Hammons} \emph{et al.}, a new era has begun. In that paper it is proved that all (non--linear) binary Kerdock-, Preparata-, Goethals- and Delsarte-Goethals-codes are images of $\mathbb{Z}_4$--linear codes under the Gray map. This discovery was quite peculiar and the paper got in 1995 the Information Theory Paper Award from the IEEE Information Theory Society. It was the start of a search for new codes by considering linear codes over the ring $\mathbb{Z}_4$ and over more \mbox{general} finite rings (see e.g.~\cite{Bini}). Also some papers by Aleksandr {\sc Nechaev} \cite{Nechaev1,Nechaev2}  have led the research in that direction.\\
We will not give a survey of all results obtained up to now for codes over finite rings, because even this niche has become too wide. We restrict ourselves here (and in the literature list) to the publications in which the direct relation between codes over rings and ring geometries is exhibited. Indeed, linear codes over finite chain rings can be associated with finite projective Hjelmslev geometries, the hyperplanes corresponding to the codewords. This correspondence offers opportunities for investigating the structure and the construction of ring--linear codes by pure geometrical methods.\\ The technique was first applied in \cite{HL1} by Thomas {\sc Honold}, now working at the Zhejiang Universit\"at in Hangzhou (China), and Ivan {\sc Landjev} from the New Bulgarian University of Sofia (Bulgaria). They prove that certain \mbox{MacDonald} codes can be represented by linear codes over the ring of twisted dual numbers on a finite field, using multisets of points in Hjelmslev spaces. 
In \cite{HL2} they prove that all Reed--Muller codes are linearly representable over the ring of dual numbers over $\mathbb{Z}_2$.
In \cite{HL3} a general theory of linear codes over finite chain rings has been developed as a natural generalization of the theory of linear codes over finite fields and the correspondence with Hjelmslev spaces is investigated. In \cite{HL4} and  \cite{HL5} an update of that paper is given. Geometric arguments are also used explicitly in \cite{HL6} for the construction of particular linear codes over chain rings of order four, genera\-lizing a result obtained by Michael {\sc Kiermaier} and Johannes {\sc Zwanzger} in \cite{K1,K2}. Keisuke {\sc Shiromoto} and Leo {\sc Storme} defined in \cite{ShSt} a Griesmer type bound for linear codes over finite quasi-Frobenius rings and they give a geometrical characterization of linear codes meeting the bound, viz.~a one-to-one correspondence between these codes and minihypers in projective Hjelmslev spaces.\\
Kiermaier was a student of Honold and  wrote his Master's thesis on {\sl ``Arcs und Codes \"uber endlichen Kettenringen"} in 2006 and obtained his PhD in 2012 at the University of Bayreuth (Germany) with the thesis {\sl ``Geometrische Konstruktionen linearer Codes \"uber Galois-Ringen der Charakteristik 4 von hoher homogener Minimaldistanz"}. 

The study of codes over chain rings from the viewpoint of Hjelmslev geometry has also led to the generalization of several special point sets (arcs, ovals, blocking sets, caps), already well--known in classical geometry over fields.  The thus far most studied sets in Hjelmslev planes are arcs. A $(k, n)$-arc refers to a set  of $k$ points which meet every line in at most $n$ points. This definition was given by Honold and Kiermaier in \cite{HL7}. General upper bounds on the cardinality of such arcs were found as well as the maximum possible size. For a chain ring $R$ of length 2 with Jacobson radical $R_0$, such that $|R/R_0|= q$, the maximum size of a $(k,2)$-arc in the projective Hjelmslev plane over $R$ is $q^2$ if $q$ is odd and $q^2 + q + 1$ if $q$ is even (see \cite{HL8}). In \cite{HL9} the existence of maximal $(q^2 + q + 1,2)$-arcs (i.e. hyperovals) is proved for $q$ even and in \cite{HK1} the existence of maximal $(q^2,2)$-arcs for $q$ odd is proved. The results on maximal arcs are also used to construct interesting codes with a linear representation over a chain ring. Examples, non--existence results and upper bounds for the length of arcs are also present in \cite{K4,HHL,K5,K7,K3,HL11,K6,KK1,K,SS}.\\
Caps in finite projective Hjelmslev spaces over chain rings of nilpotency index 2 are defined by Honold and Landjev in \cite{HL10}. A geometric construction for caps in the threedimensional space is given, using ovoids in the epimorphic space PG(3,$q$) as well as an algebraic construction using Teichm\"uller sets. Blocking sets in Hjelmslev planes and their relation with codes is the subject of \cite{L4,L2,L3} while \cite{L1} and \cite{LG} deal with spreads. 

Another aspect that appears in the literature is the correspondence between two--weight codes and strongly regular graphs. In \cite{G1}  regular projective two-weight codes over finite Frobenius rings are introduced and it is shown that such codes give rise to a strongly regular graph. In \cite{G1,G2} two-weight codes are constructed using ring geometries and they yield infinite families of strongly regular graphs with non-trivial parameters.

Ovals in an ordinary projective plane of order $q$ are just $q+1$ arcs and every conic is an oval. By a celebrated theorem of Segre every $(q+1)$--arc in PG(2,$q$) with $q$ odd, is a conic.
The study of ovals, conics and unitals in finite projective Klingenberg and Hjelmslev planes was initiated by the author in his PhD thesis but only a part of it was published, e.g.~in \cite{KepMi}.
The author also defined polarities in projective Klingenberg planes and spaces and he investigated their sets of absolute points which give rise to ovals or ovoids in some cases. A comprehensive study of polarities in $n$--uniform Hjelmslev planes and spaces over the ring GF($q$)[$t$]/$t^n$ appeared in \cite{KepP1,KepP2,KepP3}. Also the papers \cite{Kossel,Kulkarni2} enlight some aspects of ovals and conics in ring planes.\\
A combinatorial study of conics in finite desarguesian  Hjelmslev planes was made by Ratislav {\sc Jurga} and Viliam {\sc Chv\'al}. Their papers contain formulas for the number of interior and exterior points, tangents and secants of a conic \cite{J4,J1,J2,J3}.
A related problem is the projective equivalence of quadrics in projective Klingenberg spaces. This was first formulated in \cite{Egory} and analyzed in detail by O.A. {\sc Starikova} et al.~in \cite{Jukl1,Star4,Star1,Star2,Star3}.\\
Recently, a special class of codes (toric codes) is found to be related to affine ring geometries (Leissner planes) by {\sc Little} in \cite{Little1,Little2}, revealing still another correspondence between ring geometry and coding theory.

\section{Ring geometry in quantum information theory: Saniga and Planat}

Besides the fact that ring  geometries play a substantial role in coding theory, there is another domain in which they arise unexpectedly, namely in quantum information theory. In this important branch of quantum physics it is studied how information can be stored and retrieved from a quantum mechanical system.
In 2006 Metod {\sc Saniga} from the Astronomical Institute at the Slovak Academy of Sciences and Michel {\sc Planat} from the Universit\'e de Franche--Comt\'e (France) discovered a connection between finite ring geometry and quantum information theory \cite{SP1}. It is not clear yet whether or not the correspondence goes further than just equality in numbers of objects, but anyhow it is remarkable that ring geometries seem to play a role in quantum physics. 

The notion of mutually unbiased bases (MUBs) has turned into a cornerstone of the modern theory. Saniga and Planat observed that the basic combinatorial properties of a complete set of MUBs  of a $q$--dimensional Hilbert space $\mathcal{H}_q$ with $q=p^r$, $p$ being a prime and $r$ a positive integer, are qualitatively mimicked by the configuration of points lying on a proper conic in a projective Hjelmslev plane defined over a Galois ring of characteristic $p^2$ and rank $r$. The $q$ vectors of a basis of $\mathcal{H}_q$ correspond to the $q$ points in a neighbour class and the $q+1$ MUBs answer to the total number of pairwise disjoint neighbour classes on the conic. 
In a series of subsequently published papers, other combinatorial correspondences between concepts from quantum theory and ring geometries over finite rings (in particular projective ring lines) are observed. One of these similarities concerns the structure of the generalized Pauli group associated with a specific single $d$-dimensional qudit (qudits are generalizations tot $d$--level quantum systems of qubits which are the basic units in quantum information systems of level two). See \cite{SP11,SP9,SP14,SP16,SP3,SP8,SP5,SP6,SP7,SP4,SP10,SP13,SP12,SP17} and the survey papers \cite{SP15,SP2} for an overview on this topic.
The study of the relation with finite geometry (not only restricted to ring geometry but also in connection with small generalized polygons and polar spaces) is still ongoing, see e.g.~\cite{SP18}.

\newpage
\section{Literature on ring geometry and geometry over rings}

In the separate  bibliographic list we only mention publications in scientific journals (or books). Other references, e.g.~Master's and PhD theses which are quoted explicitly in the foregoing text, are not repeated here, except in case they were also published. Articles belonging to conference proceedings are also omitted if there exist copies with almost identical content, published in other journals. We have grouped the references in themes, corresponding to the sections in the paper. Some papers can be classified in multiple sections. In that case they are listed in the section in which they are referred to for the first time. 
We do not claim that the bibliography is complete but we have tried to compose an accurate list which complements the existing obsolete lists, especially for the period after 1990.
The author will be grateful for additions, completions and corrections to this list.\\

\renewcommand{\refname}{A \, First traces and pioneers of ring geometry}

\renewcommand{\refname}{B \, The Belgian contribution to early ring geometry}

\renewcommand{\refname}{C \, The foundations of plane affine ring geometry}

\renewcommand{\refname}{D \, Metric geometry over rings}

\renewcommand{\refname}{E \,The florescence of Hjelmslev geometry}

\renewcommand{\refname}{F \, The continuation of the Hjelmslev epoch}

\renewcommand{\refname}{G  \, The coordinatization of Klingenberg and Hjelmslev planes}

\renewcommand{\refname}{H \, Order and topology in ring geometry}

\renewcommand{\refname}{I \, The revival of ring geometry (Veldkamp, Faulkner)}

\renewcommand{\refname}{J \, Projective and affine Hjelmslev spaces and spaces over rings}

\renewcommand{\refname}{K \,Projective lines and circle geometries over rings}

\renewcommand{\refname}{L \, Ring geometries and buildings}

\renewcommand{\refname}{M \,Ring geometries and coding theory}

\renewcommand{\refname}{N \,Ring geometries in quantum information theory}

\end{document}